\documentclass[12pt]{amsart}
\usepackage{amscd}
\usepackage[all]{xy}

\def\NZQ{\Bbb}               
\def\NN{{\NZQ N}}

\def\CC{{\NZQ C}}

\def\frk{\frak}               

\def\mm{{\frk m}}

\def\Phi{{\frk n}}
\def\Phi{{\frk N}}

\def\opn#1#2{\def#1{\operatorname{#2}}} 
\opn\chara{char} \opn\length{\ell} \opn\pd{pd} \opn\rk{rk}
\opn\projdim{proj\,dim} \opn\injdim{inj\,dim} \opn\rank{rank}
\opn\depth{depth} \opn\grade{grade} \opn\height{height}
\opn\embdim{emb\,dim} \opn\codim{codim}

\opn\Tr{Tr} \opn\bigrank{big\,rank}
\opn\superheight{superheight}\opn\lcm{lcm}
\opn\trdeg{tr\,deg}
\opn\reg{reg} \opn\lreg{lreg} \opn\ini{in} \opn\lpd{lpd}
\opn\size{size} \opn\Pf{Pf} \opn\GL{GL} \opn\SL{SL} \opn\mod{mod}
\opn\ord{ord} \opn\Gin{Gin}
\opn\Hilb{Hilb}\opn\adeg{adeg}\opn\std{std}\opn\ip{infpt}
\opn\pol{pol}
\opn\div{div} \opn\Div{Div} \opn\cl{cl} \opn\Cl{Cl}
\opn\Spec{Spec} \opn\Supp{Supp} \opn\supp{supp} \opn\Sing{Sing}
\opn\Ass{Ass} \opn\Min{Min}
\opn\Ann{Ann} \opn\Rad{Rad} \opn\Soc{Soc}
\opn\Syz{Syz} \opn\Im{Im} \opn\Ker{Ker} \opn\Coker{Coker}
\opn\Am{Am} \opn\Hom{Hom} \opn\Tor{Tor} \opn\Ext{Ext}
\opn\End{End} \opn\Aut{Aut} \opn\id{id}

\opn\nat{nat}
\opn\pff{pf}
\opn\Pf{Pf} \opn\GL{GL} \opn\SL{SL} \opn\mod{mod} \opn\ord{ord}
\opn\Gin{Gin} \opn\Hilb{Hilb}
\opn\and{and}
\opn\aff{aff} \opn\con{conv} \opn\relint{relint} \opn\st{st}
\opn\lk{lk} \opn\cn{cn} \opn\core{core} \opn\vol{vol}
\opn\link{link} \opn\star{star}
\opn\gr{gr}

\def\pot#1#2{#1[\kern-0.28ex[#2]\kern-0.28ex]}

\opn\dirlim{\underrightarrow{\lim}}
\opn\inivlim{\underleftarrow{\lim}}
\let\to=\rightarrow

\def\Implies{\ifmmode\Longrightarrow \else
        \unskip ${}\Longrightarrow{} $\ignorespaces\fi}
\def\implies{\ifmmode\Rightarrow \else
        \unskip ${}\Rightarrow{}$\ignorespaces\fi}
\def\iff{\ifmmode\Longleftrightarrow \else
        \unskip ${}\Longleftrightarrow{}$\ignorespaces\fi}

\let\:=\colon
\newtheorem{Theorem}{Theorem}[section]
\newtheorem{Lemma}[Theorem]{Lemma}
\newtheorem{Corollary}[Theorem]{Corollary}
\newtheorem{Proposition}[Theorem]{Proposition}
\newtheorem{Remark}[Theorem]{Remark}

\theoremstyle{definition}
\newtheorem{Example}[Theorem]{Example}

\let\epsilon\varepsilon
\let\phi=\varphi
\let\kappa=\varkappa
\def\qed{\ifhmode\textqed\fi
      \ifmmode\ifinner\quad\qedsymbol\else\dispqed\fi\fi}
\def\textqed{\unskip\nobreak\penalty50
       \hskip2em\hbox{}\nobreak\hfil\qedsymbol
       \parfillskip=0pt \finalhyphendemerits=0}
\def\dispqed{\rlap{\qquad\qedsymbol}}

\opn\dis{dis}
\def\pnt{{\raise0.5mm\hbox{\large\bf.}}}

\opn\Lex{Lex}

\newcommand{\gin}{\ensuremath{\mathrm{Gin}}}

\newcommand{\init}{\ensuremath{\mathrm{in}}}
\newcommand{\Image}{\ensuremath{\mathrm{Im}}}
\newcommand{\mideal}{\ensuremath{\mathfrak{m}}}

\def\cocoa{{\hbox{\rm C\kern-.13em o\kern-.07em C\kern-.13em o\kern-.15em A}}}

\newcommand{\M}{ \ensuremath{\mathcal{M}}}

\begin{document}

\title
{Rigidity of Linear Strands and Generic Initial Ideals}

\author{Satoshi Murai}
\address{
Department of Pure and Applied Mathematics\\
Graduate School of Information Science and Technology\\
Osaka University\\
Toyonaka, Osaka, 560-0043, Japan\\
}
\email{s-murai@ist.osaka-u.ac.jp
}
\author{Pooja Singla}
\address{FB Mathematik, Universit\"at Duisburg-Essen, Campus Essen, 45117 Essen, Germany}
\email{pooja.singla@uni-duisburg-essen.de}

\thanks{The first author is supported by JSPS Research Fellowships for Young Scientists}

\maketitle

\begin{abstract}
Let $K$ be a field, $S$ a polynomial ring and $E$ an exterior
algebra over $K$, both in a finite set of variables. We study rigidity properties of
the graded Betti numbers of graded ideals in $S$ and $E$ when passing to their generic initial ideals. 
First, we prove that if
the graded Betti numbers $\beta_{ii+k}^S(S/I)=\beta_{ii+k}^S\big(S/\Gin(I)\big)$ for some $i>1$ 
and $k \geq 0$, then  $\beta_{qq+k}^S(S/I)= \beta_{qq+k}^S\big(S/\Gin(I)\big)$ for all $q \geq i$,
where $I\subset S$ is a graded ideal. Second, we show that if
$\beta_{ii+k}^E(E/I)= \beta_{ii+k}^E\big(E/\Gin(I)\big)$ for some $i>1$ and $k \geq 0$,
then  $\beta_{qq+k}^E(E/I)= \beta_{qq+k}^E\big(E/\Gin(I)\big)$ for all $q \geq 1$, where
$I\subset E$ is a graded ideal. In addition, it will be shown that the graded Betti numbers
$\beta_{ii+k}^R(R/I)= \beta_{ii+k}^R\big(R/\Gin(I)\big)$ for all $i \geq 1$
if and only if $I_{\langle k \rangle}$ and $I_{\langle k+1 \rangle}$ have a linear resolution. Here
$I_{\langle d \rangle}$ is the ideal generated by all homogeneous elements in $I$ of degree $d$,
and $R$ can be either the polynomial ring or the exterior algebra.
\end{abstract}
\section*{Introduction}
In this paper we study rigidity properties of graded Betti numbers of a graded ideal when passing 
to its generic initial ideal.

Let $S=K[x_1,\dots,x_n]$ be the polynomial ring in $n$ variables over a field $K$ with char$(K)=0$ and
$I \subset S$ a graded ideal. Let $\beta_{i}^S(M)=\dim_K\Tor_{i}^S(K,M)$ and 
$\beta_{ij}^S(M)=\dim_K\Tor_{i}^S(K,M)_j$ denote respectively the $i$-th total and
{$i,j$-th graded Betti number}  of a finitely generated graded $S$-module $M$.\pagebreak[3]

The generic initial ideal $\Gin(I)$ plays a fundamental role in investigating
various homological, algebraic, combinatorial and geometric properties of $I$.
By definition, the generic initial ideal $\Gin(I)$ is, after performing a generic change of coordinates, the initial ideal of $I$
with respect to the reverse lexicographic order.
Here we consider the reverse lexicographic order
induced by $x_1> \cdots >x_n$.\pagebreak[3]

The following inequality of the graded Betti numbers is well-known:
$$\beta_{ij}(S/I)\leq \beta_{ij}\big(S/\Gin(I)\big),$$
for all $i,j$ (see \cite[Theorem 1.1]{C}). Equality holds for all
$i$ and $j$ if and only if $I$ is componentwise linear (see \cite[Theorem 1.1]{Stable}).
In his paper \cite{C} Conca  asked whether the equality $\beta_{i}(S/I) = \beta_{i}\big(S/\Gin(I)\big)$ 
for some $i \geq 1$ of the total Betti numbers implies
$\beta_{j}(S/I)= \beta_{j}\big(S/\Gin(I)\big)$ for all $j \geq i$. This question of Conca
was positively answered in 2004 by Conca, Herzog and Hibi in \cite{CHH}.\pagebreak[3]

One of the main results of our paper is to extend this result of  Conca--Herzog--Hibi to graded Betti numbers.
In Corollary \ref{rigidgin}
we show the following: If for some $i > 1$ and  $k \geq 0$, we have
$\beta_{ii+k}^S(S/I)=\beta_{ii+k}^S\big(S/\Gin(I)\big)$,
then
$$\beta_{qq+k}^S(S/I)=\beta_{qq+k}^S\big(S/\Gin(I)\big)\ \ \  \mbox{ for all } q \geq i.$$ 

We also study the same property for generic initial ideals over an exterior algebra.
Let $K$ be an infinite field,
$V$ an $n$-dimensional $K$-vector space with basis $e_1,\dots,e_n$
and $E= \bigoplus_{k=0}^{n} \bigwedge^k V$ the exterior algebra of $V$.
For a graded ideal $J\subset E$,
we write $\Gin(J)$ for the generic initial ideal of $J$ with respect to the reverse lexicographic order
induced by $e_1>\dots >e_n$
and denote by $\beta_{ij}^E(E/J)$
the {$i,j$-th graded Betti number} of $E/J$ over $E$.
Somewhat surprisingly,
the following stronger property is true in the exterior algebra:
If $\beta_{ii+k}^E(E/J) = \beta_{ii+k}^E \big(E/\Gin(J)\big)$ for some $i > 1$ and  $k \geq 0$,
then one has
$$\beta_{qq+k}^E(E/J) = \beta_{qq+k}^E \big(E/\Gin(J)\big) \ \ \ \ \ \mbox{ for all }\ \ \ 
q \geq 1.$$

Let $R$ be either a polynomial ring over a field $K$ with $\mathrm{char}(K)=0$
or an exterior algebra over an infinite field
and $I$ a graded ideal of $R$.
The above property leads us to ask when a graded ideal $I \subset R$
satisfies $\beta_{ii+k}^R(R/I) = \beta_{ii+k}^R \big(R/\Gin(I)\big)$ for all $i \geq 1$,
where we fix an integer $k \geq 0$.
We will prove the following result answering this question.\pagebreak[3]

\begin{Theorem} \label{dlinear}
Let $R$ be either a polynomial ring over a field $K$ with $\mathrm{char}(K)=0$
or an exterior algebra over an infinite field,
$I \subset R$ a graded ideal and $k \geq 0$ an integer.
The following conditions are equivalent.
\begin{itemize}
\item[(i)] $\beta_{ii+k}^R(R/I)= \beta_{ii+k}^R\big(R/\Gin(I)\big)$ for all $i\geq 1$;
\item[(ii)] $I_{\langle k \rangle}$ and $I_{\langle k+1 \rangle}$ have a linear resolution;
\item[(iii)] $\beta_{1k+1}^R(R/I)=\beta_{1k+1}^R\big(R/\Gin(I)\big)$ and 
$\beta_{1k+2}^R(R/I)=\beta_{1k+2}^R\big(R/\Gin(I)\big)$,
\end{itemize}
where $I_{\langle k \rangle}$ denotes the ideal of $R$
generated by all homogeneous elements in $I$ of degree $k$.
\end{Theorem}\pagebreak[3]

The above result is a generalization of \cite[Theorem 1.1]{Stable},
where it was shown that $\beta_{ij}^R(R/I)= \beta_{ij}^R\big(R/\Gin(I)\big)$ for all $i,j$
if and only if $I$ is componentwise linear.\pagebreak[3]

In the end of the paper, we study the Cancellation Principle for generic initial ideals \cite{G}.
We find the relation between our results for Betti numbers of a graded ideal
in a polynomial ring and the Cancellation Principle for generic initial ideals.\pagebreak[3]

This paper is organized as follows:
In \S 1,
we will give an upper bound for graded Betti numbers in terms of generic annihilator numbers
by using the technique developed in \cite{CHH}.
In \S 2,
we will generalize Conca--Herzog--Hibi's theorem for graded Betti numbers over a polynomial ring.
In \S 3,
some basic facts about Cartan complexes and generic annihilator numbers over an exterior algebra are studied.
In \S 4,
we will generalize Conca--Herzog--Hibi's theorem for graded Betti numbers over an exterior algebra.
In \S 5,
we will study when $I_{\langle d \rangle}$ has a linear resolution from the viewpoint of generic initial ideals
and give a proof of Theorem \ref{dlinear}.
In \S 6, we will study the Cancellation Principle. The results in the last section are closely related
to the results in \S 1.\pagebreak[3]

\section{An upper bound for the graded Betti numbers }
In this section,
we will give an upper bound for graded Betti numbers
in terms of generic graded annihilator numbers,
which were introduced in \cite{CHH}.
Note that most of the results in this section
are refinements of the results in \cite[\S 1]{CHH}.
Though these results seem to be somewhat technical, they are of crucial
importance for the proof of one of our main theorems in the next section.\pagebreak[3]

Let $S=K[x_1,\ldots,x_n]$ be the standard graded  polynomial ring
over an arbitrary field $K$ and $\mm=(x_1,\ldots,x_n)$ the
graded maximal ideal. Let $M$ be a finitely generated graded $S$-module.
For each nonnegative integer $i$, the modules
$\Tor_i^S(K,M)$ are finitely generated  $K$-vector spaces. The
numbers $\beta^S_{i}(M)=\dim_K\Tor_i^S(K,M)$ and
$\beta^S_{ij}(M)=\dim_K\Tor_i^S(K,M)_j$ are called {\em Betti
numbers} and {\em graded Betti numbers} of $M$, respectively.
As $\beta^S_{ij}$ are invariants under base field extensions, from
now on we may assume the field $K$ to be infinite.\pagebreak[3]

Let $y_1,\ldots,y_n$ be a sequence of generic linear forms for the module $M$.
For each $p=1,\ldots,n$, the modules
$$A_p = (y_1,\ldots,y_{p-1})M:_M y_p/(y_1,\ldots,y_{p-1})$$ are $\NN$-graded $S$-modules of finite length.
We define $\alpha_p(M)=\dim_KA_p$, which we call the {\em generic
annihilator numbers} of $M$. We denote by $\alpha_{p,j}(M)$ the vector
space dimension of the $j$th graded component $(A_{p})_j$ of $A_p$
which we call the {\em generic graded annihilator numbers} of
$M$. \pagebreak[3]

Let $H_i(p,M)$ be the Koszul homology
$H_i(y_1,\ldots,y_p;M)$ of the partial sequence $y_1,\ldots,y_p$.
We set $h_i(p,M)=\dim_K H_i(p,M)$ and $h_{ij}(p,M)=\dim_K
H_i(p,M)_j$. We omit $M$  and simply write $\beta^S_{ij}$,
$\beta^S_i$, $\alpha_{i,j}$, $\alpha_{i}$, $H_i(p)_j$, $H_i(p)$,
$h_{ij}(p)$, $h_i(p)$ for the above defined terms, if the module
under consideration is fixed.
Then we have the following long exact sequence (see \cite[Corollary 1.6.13]{BH}):
\begin{eqnarray} \label{Kos}&&
\begin{array}{cccccccc}
\cdots&\longrightarrow& H_i(p-1)&\xrightarrow{\phi_{i,p-1}} H_i(p-1)&\longrightarrow& 
H_i(p)&\longrightarrow& H_{i-1}(p-1)\\
\cdots&\longrightarrow& H_0(p-1)&\xrightarrow
{\phi_{0,p-1}} H_0(p-1)&\longrightarrow& H_0(p)&\longrightarrow& 0.\\
\end{array}
\end{eqnarray}\pagebreak[3]
In the above sequence $\phi_{i,p-1}$ is the multiplication map on
$H_i(p-1)$ with multiplication by $\pm y_p$. One may notice that $A_p$ is given by
the kernel of the map $\phi_{0,p-1}$. Hence we get the following exact sequences with all the maps of  degree zero:
$$0 \longrightarrow \Im\phi_{1,p-1}\longrightarrow H_1(p-1)\longrightarrow H_1(p) \longrightarrow A_p(-1)
\longrightarrow 0 $$
for all $p$, and
$$0 \longrightarrow \Im\phi_{i,p-1}\longrightarrow H_i(p-1)\longrightarrow H_i(p)\longrightarrow
H_{i-1}(p-1)(-1)\longrightarrow\Im\phi_{i-1,p-1}\longrightarrow 0,$$
for all $p$ and $i>1.$ \pagebreak[3]

Let $\delta_{i,j,k}=\dim_K(\Im\phi_{i,j})_{k}.$ From the above exact sequences, we obtain  the
following equations for each integer $k \geq 0$:
\begin{equation}\label{first}
h_{1k}(p)=h_{1k}(p-1)+\alpha_{p,k-1}-\delta_{1,p-1,k},
\end{equation}\pagebreak[3]
and  for all $i>1$,
\begin{equation}\label{second}
h_{i,i+k}(p)=h_{i,i+k}(p-1)+h_{i-1,i-1+k}(p-1)-
\delta_{i,p-1,i+k}-\delta_{i-1,p-1,i+k}.
\end{equation}\pagebreak[3]
By using (\ref{first}) and (\ref{second}), we obtain

\begin{Proposition}\label{graded}
 For all  nonnegative integers $i \geq 1$ and $k$, one has
 \begin{eqnarray}
h_{i,i+k}(p)&=&\sum_{j=1}^{p-i+1}\binom{p-j}{i-1}\alpha_{j,k}\\
\nonumber&&-\sum_{(a,b)\in A_{i,p}}\left[\binom{p-b-1}{i-a}\delta_{a,b,a+k}
+\binom{p-b-1}{i-a-1}\delta_{a,b,a+k+1}\right],
\end{eqnarray}
 where the set $A_{i,p}=\big\{(a,b)\in\NN^2 : 1\leq b\leq p-1 \;\and \; \max\{i-p+b,1\}\leq a \leq i\big\}$.
\end{Proposition}
\begin{proof} We will prove the above formula by induction on $p$.
For $p=1$, we have from Equation (\ref{first}) and Equation (\ref{second}):
$$h_{i,i+k}(1) = \begin{cases}
{\alpha_{1,k}} & {\mbox{if} \;i=1},\\
{0}&{i\geq 2}.\end{cases} $$
which is what the formula given in the statement of the proposition suggests.
Now we assume $p>1$ and we assume the result to be true for $p-1$.\pagebreak[3]

Let first $i=1.$ By induction hypothesis and from
Equation (\ref{first}), we get :
\begin{eqnarray}
h_{1,1+k}(p) \!\! &=& \! h_{1,1+k}(p-1)+\alpha_{p,k}-\delta_{1,p-1,1+k}\nonumber \\
&=& \! \sum_{j=1}^{p-1}\binom{p-1-j}{0}\alpha_{j,k}- \sum_{(a,b)\in
A_{1,p-1}}\binom{p-b-2}{1-a}\delta_{a,b,a+k}
+\alpha_{p,k}-\delta_{1,p-1,1+k}\nonumber \\
&=& \! \sum_{j=1}^{p}\alpha_{j,k}-\sum_{(a,b)\in
A_{1,p}}\left[\binom{p-b-1}{1-a} \delta_{a,b,a+k}\right] \label{AC3}\nonumber
\end{eqnarray}
which is what the formula suggests.\pagebreak[3]

Now let $i>1$. From
Equation (\ref{second}), we have:
\begin{eqnarray*}
h_{i,i+k}(p)&=&h_{i,i+k}(p-1)+h_{i-1,i-1+k}(p-1)-\delta_{i,p-1,i+k}-\delta_{i-1,p-1,i+k}.
\end{eqnarray*}
Note that one has ${a \choose b} + {a \choose b+1} ={a+1 \choose b+1}$
for all integers $a \geq b \geq 0$.
Then,
using induction hypothesis,
the right hand side of the above equation is a sum of the following three terms:
\begin{eqnarray}\label{third}
& \sum_{j=1}^{p-i+1}\big\{ { p-j-1 \choose i-1} + {p-j-1 \choose i-2}\big\} \alpha_{j,k}
=   \sum_{j=1}^{p-i+1}\binom{p-j}{i-1} \alpha_{j,k},
\end{eqnarray}\pagebreak[3]
\begin{eqnarray}\label{fourth}
&-\big\{\sum_{(a,b)\in A_{i,p-1}}{\binom{p-b-2}{i-a}}\delta_{a,b,a+k}+\delta_{i,p-1,i+k}
+\sum_{(a,b)\in A_{i-1,p-1}}\binom{p-b-2}{i-a-1}\delta_{a,b,a+k}\big\},
\end{eqnarray}\pagebreak[3]
and
\begin{eqnarray}\label{fifth}
&-\big\{\sum_{(a,b)\in A_{i,p-1}}\binom{p-b-2}{i-a-1}\delta_{a,b,a+k+1}+\delta_{i-1,p-1,i+k}\\
&+\sum_{(a,b)\in A_{i-1,p-1}}\binom{p-b-2}{i-a-2}\delta_{a,b,a+k+1}\big\}.\nonumber
\end{eqnarray}\pagebreak[3]
The term (\ref{fourth}) can be written as:

\begin{center}
$-\big\{\sum_{(a,b)\in A_{i,p-1}}\binom{p-b-2}{i-a}\delta_{a,b,a+k}
+\sum_{(a,b)\in A_{i,p-1}}\binom{p-b-2}{i-a-1}\delta_{a,b,a+k}
+\delta_{i ,p-1,i+k}\big\},$
\end{center}\pagebreak[3]
which is further equal to
 \begin{center}
$-\big\{\sum_{(a,b)\in A_{i,p-1}}\binom{p-b-1}{i-a}\delta_{a,b,a+k}+\delta_{i,p-1,i+k}\big\},$
\end{center}\pagebreak[3]
which in the end equals 
\begin{eqnarray}\label{sixth}
&-\sum_{(a,b)\in A_{i,p}}\binom{p-b-1}{i-a}\delta_{a,b,a+k}.
\end{eqnarray}\pagebreak[3]
Now we notice that the term (\ref{fifth}) can be written as:
\begin{eqnarray*}
&-\big\{\sum_{(a,b)\in
A_{i,p-1}}\binom{p-b-2}{i-a-1}\delta_{a,b,a+k+1}
+\sum_{(a,b)\in A_{i,p-1}}\binom{p-b-2}{i-a-2}\delta_{a,b,a+k+1}\\
&+\sum_{b=p-i+1}^{p-2}\delta_{i-p+b,b,i-p+b+k+1}+
\delta_{i-1,p-1,i+k}\big\}.
\end{eqnarray*}\pagebreak[3]
This can be rewritten as:
\begin{eqnarray*}
&-\big\{\sum_{(a,b)\in A_{i,p-1}}\binom{p-b-1}{i-a-1}
\delta_{a,b,a+k+1}
+\sum_{b=p-i+1}^{p-2}\delta_{i-p+b,b,i-p+b+k+1}+\delta_{i-1,p-1,i+k}\big\},
\end{eqnarray*}\pagebreak[3]
which then is equal to
\begin{eqnarray}\label{seventh}
&-\sum_{(a,b)\in A_{i,p}}\binom{p-b-1}{i-a-1}\delta_{a,b,a+k+1}.
\end{eqnarray}\pagebreak[3]
Hence
$h_{i,i+k}(p)$ is the sum of (\ref{third}),(\ref{sixth}) and (\ref{seventh}), as required.
\end{proof}
\begin{Remark}{\em  Notice that summing the formula
stated in Proposition \ref{graded} over $k$, gives us back the
formula given in the proof of \cite[Proposition 1.1]{CHH}}.
\end{Remark}\pagebreak[3]
Proposition \ref{graded} implies the following fact.
\begin{Corollary}\label{h}
We have
\begin{enumerate}
\item[(a)] $h_{i,i+k}(p) \leq \sum_{j=1}^{p-i+1}\binom{p-j}{i-1}\alpha_{j,k}$.
\item[(b)] For given integers $i \geq 1$ and $p\geq 1$, the following conditions are equivalent:
\begin{enumerate}
\item[(i)] $h_{i,i+k}(p)=\sum_{j=1}^{p-i+1}\binom{p-j}{i-1}\alpha_{j,k}$
\item[(ii)] $(\Im\phi_{a,b})_{(a+k)\phantom{+1}}=0$ for all $(a,b) \in A_{i,p}\setminus 
\big\{(i-p+b,b): b \leq p-1\big\}$
and $(\Im\phi_{a,b})_{(a+k+1)}=0$ for all $(a,b) \in A_{i,p}\setminus\big\{(i,b):
b \leq p-1\big\}$.
\item[(iii)]$\big( \mm H_a(b)\big)_{(a+k)\phantom{+1}}=0$ for all $(a,b) \in A_{i,p}\setminus
 \big\{(i-p+b,b): b \leq p-1\big\}$
and $\big(\mm H_a(b)\big)_{(a+k+1)}=0$ for all $(a,b) \in A_{i,p}\setminus\big\{(i,b):
b \leq p-1\big\}$.
\end{enumerate}
\end{enumerate}
\end{Corollary}
\begin{proof}
Statement (a) is clear from Proposition \ref{graded}. The
equivalence of (i) and (ii) follows  immediately from Proposition
\ref{graded}. Indeed,
$h_{i,i+k}(p)=\sum_{j=1}^{p-i+1}\binom{p-j}{i- 1}\alpha_{j,k}$ if
and only if all  graded maps appearing in the formula in
Proposition \ref{graded} vanish whenever their  binomial
coefficients are nonzero. And for the equivalence of (ii) and
(iii), we may notice that a generic linear form annihilates
$\big(H_a(b)\big)_k$ if and only if $\mm$ annihilates  $\big(H_a(b)\big)_{k}$.
 \end{proof}

The next corollary is a special case $(p=n)$ of
the above corollary.\pagebreak[3]
\begin{Corollary}\label{b}
\begin{enumerate}
\item[(a)] $\beta^S_{ii+k}\leq \sum_{j=1}^{n-i+1}\binom{n-j}{i-1}\alpha_{j,k}$ for all $i\geq 1$.
\item[(b)]
For given $i \geq 1$ the following are equivalent:
\begin{enumerate}
\item[(i)]$\beta^S_{ii+k}=\sum_{j=1}^{n-i+1}\binom{n-j}{i-1}\alpha_{j,k}$.
\item[(ii)]
 $(\Im\phi_{a,b})_{(a+k)\phantom{+1}}=0$ for all $(a,b) \in A_{i,n}\setminus \big\{(i-n+b,b), b \leq n-1\big\}$
and $(\Im\phi_{a,b})_{(a+k+1)}=0$ for all $(a,b) \in A_{i,n}\setminus\big\{(i,b),
b \leq n-1\big\}$.
\item[(iii)] $\big( \mm H_a(b)\big)_{(a+k)\phantom{+1}}=0$ for all $(a,b) \in A_{i,n}\setminus 
\big\{(i-n+b,b), b \leq n-1\big\}$
and $\big(\mm H_a(b)\big)_{(a+k+1)}=0$ for all $(a,b) \in A_{i,n}\setminus\big\{(i,b), b \leq n-1\big\}.$
\end{enumerate}
\end{enumerate}
\end{Corollary}\pagebreak[3]
\section{Graded Rigidity of Resolutions and Linear Components}
In this section we generalize \cite[Theorem 2.3]{CHH} of
Conca--Herzog--Hibi. They gave an upper bound of total Betti numbers in terms of generic annihilator numbers,
and proved that if the Betti number
$\beta^S_{i}(M)$ for some $i \geq 1$ reaches its upper bound, then
the Betti numbers $\beta^S_q(M)$  also reach
their upper bounds for all $q \geq i$. We show that if a graded
Betti number $\beta^S_{ii+k}(M)$ for some $i>1$ reaches its upper
bound given in Corollary \ref{b}, then so do all the graded Betti numbers $\beta^S_{qq+k}(M)$
for $q \geq i$. Here we need the assumption $i>1$ as we will see later
in Remark \ref{why}.\pagebreak[3]

We state the main theorem of this section:
\begin{Theorem} \label{jaan} Let $M$ be a finitely generated graded $S$-module. Suppose for some $i >1$, we have
$\beta^S_{ii+k}(M)=\sum_{j=1}^{n-i+1}\binom{n-j}{i-1}\alpha_{j,k}(M)$. Then
$$\beta^S_{qq+k}(M)=\sum_{j=1}^{n-q+1}\binom {n-j}{q-1}\alpha_{j,k}(M) \; \mbox{for\;all}\; q \geq i.$$
\end{Theorem}\pagebreak[3]
Before proving the theorem,
we recall the following vanishing property
of Koszul homology.
For a sequence of elements $y_1,\ldots,y_r \in S$
and a set $A \subseteq \{1,\ldots,r\}$, we set
$y_A=\{y_j : j \in A\}$.\pagebreak[3]
\begin {Lemma}\label{gcor} Let $I \supseteq (y_1,\ldots,y_r) $ and assume that $
\big(IH_i(y_A; M)\big)_{i+k}=0$ for all $A \subseteq \{1,\ldots,r\}$ for some $i,k$. Then
$\big(IH_{i+1}(y_A,M)\big)_{i+k+1}=0$ for all $A \subseteq \{1,\ldots,r\}$.
\end{Lemma}\pagebreak[3]
The proof of Lemma \ref{gcor} is the same as \cite[Corollary 2.3]{CHH}.
Hence we skip the proof.
\begin{proof}[Proof of Theorem \ref{jaan}]
First we notice that it is  enough to prove
the claim in the case when $q=i+1$. Therefore we only need to show that $\big(
\mm H_a(b)\big)_{a+k}=0$ for all $(a,b) \in A_{i+1,n}\setminus
\big\{(i+1-n+b,b): b \leq n-1\big\}$ and $\big(\mm H_a(b)\big)_{a+k+1}=0$ for all
$(a,b) \in A_{i+1,n}\setminus\big\{(i+1,b): b \leq n-1\big\}$, as is
clear from Corollary \ref{h}.\pagebreak[3]

By assumption, $\big( \mm H_a(b)\big)_{a+k}=0$ for
all $(a,b) \in A_{i,n}\setminus \big\{(i-n+b,b): b \leq n-1\big\}$ and
$\big(\mm H_a(b)\big)_{a+k+1}=0$ for all $(a,b) \in
A_{i,n}\setminus\big\{(i,b):b \leq n-1\big\}$.
Also a routine computation implies
$$A_{i+1,n}\setminus \big(A_{i,n}\setminus \big\{(i-n+b,b): b \leq
n-1\big\}\big)=\big\{(i+1,b): b \leq n-1\big\}$$
and
$$\big(A_{i+1,n}\setminus\{(i+1,b): b \leq n-1\}\big)\setminus
\big(A_{i,n}\setminus\{(i,b): b \leq n-1\}\big)=\big\{(i,b): b \leq
n-1\big\}.$$
Therefore, we need to show that $\big(\mm
H_{i+1}(b)\big)_{i+1+k}=0$  and $\big(\mm H_i(b)\big)_{i+k+1}=0$ for all $b
\leq n-1$.
However, from assumption $\big(\mm H_{i}(b)\big)_{i+k}=0$ and $\big(\mm
H_{i-1}(b)\big)_{i+k}=0$ for all $b \leq n-1$, now it follows from Lemma
\ref{gcor} that for all $b \leq n-1$, we have $\big(\mm H_{i+1}(b)\big)_{i+1+k}=0$  and $\big(\mm
H_i(b)\big)_{i+k+1}=0$. Hence we are done.
\end{proof}\pagebreak[3]
A graded ideal $I \subset S$ generated in degree $d$ is said to have a
{\em linear resolution} if the regularity $\reg(I)=\max\big\{k:\beta^S_{ii+k}(I) \ne 0\big\}$ of $I$ is
equal to $d$.
Also, a graded ideal $I$ is said to be {\em componentwise
linear}  if the ideal $I_{\langle k\rangle}$ has linear resolution
for each $k$.
A monomial ideal $I \subset S$ is said to be \textit{strongly stable}
if $ux_q \in I$ implies $ux_p\in I$ for any $1 \leq p<q \leq n$.
Note that generic initial ideals are strongly stable if $\mathrm{char}(K) =0$,
and strongly stable ideals are componentwise linear.\pagebreak[3]

Theorem \ref{jaan} has a nice meaning in the special case
$M=S/I$ where $I$ is a graded ideal in $S$.
Let $I \subset S$ be a graded ideal and $\Gin(I)$ its generic
initial ideal with respect to the reverse lexicographic order.
It follows from \cite[Theorem 1.5]{CHH} that a graded ideal $I \subset S$
is componentwise linear if and only if the Betti numbers of $S/I$ reaches the upper bound
given in Corollary \ref{b}.
Also, it is not hard to show that $\alpha_{i,j}(S/I)=\alpha_{i,j}(S/\Gin(I))$
for all $i$ and $j$
(see \cite[Lemma 2.5]{CHH}).\pagebreak[3]
Then, since $\Gin(I)$ is componentwise linear,
we have 
\begin{eqnarray}\label{TT}
\beta^S_{ii+k}\big(S/\Gin(I)\big)=\sum_{j=1}^{n-i+1}\binom{n-j}{i-1}\alpha_{j,k}(S/I)\ \ \mbox{ for all } 
i \mbox{ and } k.
\end{eqnarray}\pagebreak[3]
This fact and Theorem \ref{jaan} immediately imply
\begin{Corollary} \label{rigidgin}
Suppose $\chara K =0$.
Let $I\subset S$ be a graded ideal.
If for some $i > 1$ and $k \geq 0$,
$\beta^S_{ii+k}(S/I)=\beta^S_{ii+k}\big(S/\Gin(I)\big)$, then
$$\beta^S_{qq+k}(S/I)=\beta^S_{qq+k}\big(S/\Gin(I)\big)  \mbox{\;for \;all \;} q \geq i.$$
\end{Corollary}
\begin{Remark} \label{why} {\em The assumption $i>1$ in Theorem \ref{jaan} (and Corollary \ref{rigidgin}) is
necessary.  In the case when $i=1$, we notice from the proof that
we need to show that $\big(\mm H_2(b)\big)_{2+k}=0$ and $\big(\mm
H_1(b)\big)_{2+k}=0$ for all $b \leq n-1$. As the set
$A_{1,n}\setminus \big\{(1,b),\; b \leq n-1\big\}= \emptyset,$ the second
equality does not follow. Moreover in the case when $M=R/I$ where
$I \subset S$ is a graded ideal, we always have
$\beta^S_{1d_0}(R/I)=\beta^S_{1d_0}\big(R/\Gin(I)\big)=\sum_{j=1}^{n+1}\alpha_{j,d_0-1}(R/I)$
where $d_0$ is the minimum of the degrees of generators of $I$.
So if Theorem \ref{jaan} would
have been true for $i=1$, then it would follow that
$\beta^S_{ii+d_0-1}(R/I)=\sum_{j=1}^{n-i+1}\binom{n-j}{i-1}\alpha_{j,d_0-1}(R/I)$
for all $i$, which is false in general.}
\end{Remark}\pagebreak[3]

As we see in Remark \ref{why},
Corollary \ref{rigidgin} is false for $i=1$.
However, the following property is true for the first graded Betti numbers.\pagebreak[3]

\begin{Corollary} \label{generator}
Suppose $\chara K =0$.
Let $I\subset S$ be a graded ideal.
Then, for a given integer $k$, the graded Betti numbers
$\beta^S_{ii+k}(S/I)= \beta^S_{ii+k}\big(S/\Gin(I)\big)$ for all $i \geq 1$
if and only if
$\beta^S_{1,k+1}(S/I)=\beta^S_{1,k+1}\big(S/\Gin(I)\big)$ and $\beta^S_{1,k+2}(S/I)=
\beta^S_{1,k+2}\big(S/\Gin(I)\big)$.
\end{Corollary}\pagebreak[3]
\begin{proof}
First, we will show the ``if'' part.
Since $\beta_{1,k+1}^S(S/I)=\beta_{1,k+1}^S\big(S/\Gin(I)\big)$
and $\beta_{1,k+2}^S(S/I)=\beta_{1,k+2}^S\big(S/\Gin(I)\big)$,
Corollary \ref{b} says that
$\mm H_1(b)_{1+k}=0$ and $\mm H_{1}(b)_{2+k}=0$ for all $b \leq n-1$.
Thus Lemma \ref{gcor} says that
$\mm H_a(b)_{a+k}=0$ and $\mm H_a(b)_{a+k+1}=0$ for all $(a,b)$ with $a \in \mathbb{Z}$ and $b \leq n-1$.
Then, by Corollary \ref{b},
we have $\beta^S_{ii+k}(S/I)= \beta^S_{ii+k}\big(S/\Gin(I)\big)$ for all $i \geq 1$.\pagebreak[3]

Next, we will show the ``only if'' part.
Since $\beta^S_{1,k+1}(S/I)=\beta^S_{1,k+1}\big(S/\Gin(I)\big)$ follows from the assumption,
what we must prove is $\beta^S_{1,k+2}(S/I)=\beta^S_{1,k+2}\big(S/\Gin(I)\big)$.
Since $\beta^S_{2,k+2}(S/I)=\beta^S_{2,k+2}\big(S/\Gin(I)\big)$,
Corollary \ref{b} says that
$\mm H_a(b)_{(a+k+1)}=0$ for all $(a,b) \in A_{2,n} \setminus \big\{ (2,b):b \leq n-1\big\}=A_{1,n}$.
This fact and Corollary \ref{b} imply $\beta^S_{1,k+2}(S/I)=\beta^S_{1,k+2}\big(S/\Gin(I)\big)$.\pagebreak[3]
\end{proof}
For any monomial $u \in S$,
write $m(u)$ for the maximal integer $i$ such that $x_i$ divides $u$. We recall a result of 
 Eliahou--Kervaire \cite{EK} which we need in the proof of our next proposition. They
proved that if $I \subset S$ is a strongly stable ideal then
\begin{eqnarray} \label{EKf}
\beta_{ii+j}(I)=\sum_{u \in G(I),\ \deg(u)=j}
{  m(u)-1 \choose i} \ \ \mbox{ for all }i \mbox{ and }j
\end{eqnarray}
where $G(I)$ is the set of minimal monomial generators of $I$.
Aramova--Herzog--Hibi \cite[Theorem 1.1]{Stable} proved that
a graded ideal $I$ in $S$ with char$(K)=0$ is componentwise linear
if and only if $\beta^S_{ij}(I)= \beta_{ij}^S(\Gin(I))$ for all $i,j$.
We will refine this result
in terms of the maximal degree of minimal generators.
\begin{Proposition}
Suppose $\chara K =0$.
Let $I\subset S$ be a graded ideal,
and let $d$ be the maximum of the degrees of the generators of $I$.
Then the following conditions are equivalent.
\begin{itemize}
\item[(i)]
$I$ is componentwise linear;
\item[(ii)]
$\beta^S_{ii+k}(I)= \beta^S_{ii+k} \big(\Gin(I)\big)$
for all $i \geq 0$ and all $k \leq d$;
\item[(iii)]
$\beta^S_{11+k}(I) = \beta^S_{11+k}\big(\Gin(I)\big)$ for all $k \leq d$;
\item[(iv)]
$\beta^S_{0k}(I)= \beta^S_{0k}\big(\Gin(I)\big)$ for all $k \leq d+1$.
\end{itemize}
\end{Proposition}
\begin{proof}
(i) $\Rightarrow$ (ii) follows from \cite[Theorem 1.1]{Stable}
and (ii) $\Rightarrow$ (iii) is obvious.
On the other hand,
 we already proved that if $\beta^S_{1k}(I)=\beta^S_{1k}\big(\Gin(I)\big)$,
 then we have $\beta^S_{0k}(I)= \beta^S_{0k}\big(\Gin(I)\big)$ in the proof of Corollary \ref{generator}.
 This fact implies (iii)$\Rightarrow$ (iv).\pagebreak[3]

Now we show (iv) $\Rightarrow$ (i). We have $\beta^S_{0d+1}(I)= \beta^S_{0d+1}\big(\Gin(I)\big)=0$, 
by assumption. Now, since $\Gin(I)$ is strongly stable, by  Eliahou--Kervaire formula (\ref{EKf}) we have
$\beta^S_{i,i+d+1}(I)=\beta^S_{i,i+d+1}\big(\Gin(I)\big)=0$ for all $i \geq 0$.\pagebreak[3]
However, the equality of graded betti numbers
$\beta^S_{1d+2}(I)=\beta^S_{1d+2}\big(\Gin(I)\big)=0$ implies the equality
$\beta^S_{0d+2}(I)= \beta^S_{0d+2}\big(\Gin(I)\big)=0$
as we see in the proof of Corollary \ref{generator}.\pagebreak[3]
Then again we have $\beta^S_{i,i+d+2}(I)=\beta^S_{i,i+d+2}\big(\Gin(I)\big)=0$ for all $i \geq 0$.
Arguing inductively,
we have $\beta^S_{0j}(I)= \beta^S_{0j}\big(\Gin(I)\big)$ for all $j \geq 0$.
Then Corollary \ref{generator} implies that $\beta_{ij}(I)= \beta_{ij}\big(\Gin(I)\big)$ for all $i,j$.
Hence
$I$ is componentwise linear.
\end{proof}\pagebreak[3]
\section{The Cartan--Complex and Generic annihilator numbers}
In this section,
we recall some basic facts about Cartan complex introduced by Cartan
and consider generic annihilator numbers in an exterior algebra.

Let $K$ be an infinite field,
$V$ an $n$-dimensional $K$-vector space with basis $e_1,\dots,e_n$
and $E=\bigoplus_{k=0}^n \bigwedge^k V$ the exterior algebra of $V$.
For any subset $S=\{i_1,\dots,i_d\}$ with $1 \leq i_1 < \cdots < i_d \leq n$,
the element $e_S=e_{i_1} \wedge \cdots \wedge e_{i_d} \in E$ is called a {\em monomial} of $E$ of degree $d$.
Let $v_1,\dots,v_m \in E_1$ be linear forms.
The \textit{Cartan complex} $C_{\bullet} (v_1,\dots,v_m;E)$ of the sequence $v_1,\dots,v_m$ is defined as the complex
whose $i$-chains $C_i (v_1,\dots,v_m;E)$ are the elements of degree $i$ of the free divided power algebra $E\langle x_1,\dots,x_m \rangle$.
Thus $C_{\bullet} (v_1,\dots,v_m;E)$ is the polynomial ring over $E$ in the set of variables
$$x_i^{(j)}, \ \quad i=1,2,\dots,m, \ \ \ j=1,2,\dots, $$
modulo the relations
$$x_i^{(j)} x_i^{(k)}= \frac{(j+k)!} {j! k!} x_i^{(k+j)},$$
where we set $x_i^{(0)}=1$ and $x_i^{(1)}=x_i$ for $i=1,\dots,m$.
The algebra $C_{\bullet} (v_1,\dots,v_m;E)$ is a free $E$-module with basis $x^{(a)}=x_1^{(a_1)}\cdots x_m^{(a_m)}$
with $a \in \mathbb{Z}^m$.\pagebreak[3]

The $E$-linear differential $\partial$ on $C_{\bullet} (v_1,\dots,v_m;E)$ is defined by
$$\partial \big(x^{(a)}\big)= \sum_{a_i>0} v_i \cdot x_1^{(a_1)} \cdots x_i^{(a_i-1)} \cdots x_m^{(a_m)}.$$
It is easily verified that $\partial \circ \partial =0$,
so that $C_{\bullet} (v_1,\dots,v_m;E)$ is indeed a complex.\pagebreak[3]

Let $\mathcal{M}$ be the category of finitely generated graded left and right $E$-module
satisfying $ax= (-1)^{\deg(a) + \deg (x)}xa$ for all homogeneous elements $a\in E$ and $x \in M$,
where $M \in \M$.
The complex $C_{\bullet} (v_1,\dots,v_m;M) =C_{\bullet} (v_1,\dots,v_m;E) \otimes_E M$ is called
the \textit{Cartan complex of $M$ with respect to $v_1,\dots,v_m \in E_1$},
and its homology $H_{\bullet} (v_1,\dots,v_m;M)$ is called the \textit{Cartan homology}.
We recall two basic properties of the Cartan homology.
(See \cite{AHH} or \cite{H} for the detail.)\pagebreak[3]
\begin{Lemma}[{\cite[Theorem 2.2]{AHH}}] \label{betti}
Let $v_1,\dots,v_n\in E$ be linearly independent linear forms
and $M \in \M$.
One has
$$ H_i(v_1,\dots,v_n;M)_j \cong \Tor^E_i(K,M)_j.$$
\end{Lemma}\pagebreak[3]
\begin{Lemma} [{\cite[Corollary 2.4]{AHH}}] \label{long}
Let $v_1,\dots,v_n\in E$ be linear forms
and $M \in \M$.
For $p=1,2,\dots,n-1$,
there exists a long exact sequence
\begin{eqnarray*}
\cdots
\stackrel{\gamma_{i,p}}{\longrightarrow}
H_i(v_1,\dots,v_p;M)
\stackrel{\eta_{i,p}}{\longrightarrow}
H_i(v_1,\dots,v_{p+1};M)
\stackrel{\psi_{i,p}}{\longrightarrow}
H_{i-1}(v_1,\dots,v_{p+1};M) (-1) \\
\stackrel{\gamma_{i-1,p}}{\longrightarrow}
H_{i-1}(v_1,\dots,v_p;M)
\stackrel{\eta_{i-1,p}}{\longrightarrow}
H_{i-1}(v_1,\dots,v_{p+1};M)
\longrightarrow
\cdots
\end{eqnarray*}
where $\eta_{i,p}$ is the map induced by the inclusion map and the maps
$\psi_{i,p}$ and $\gamma_{i,p}$ are defined as follows:
If $z=g_0+g_1x_{p+1}+\cdots+ g_ix_{p+1}^{(i)}$ is a cycle in $C_i(v_1,\dots,v_{p+1};M)$
with each $g_k \in C_i(v_1,\dots,v_{p};M)$,
then $\psi_{i,p}\big([z]\big)
=[g_1+g_2x_{p+1}+\cdots+ g_ix_{p+1}^{(i-1)}]$
and $\gamma_{i,p}\big([z]\big) =[g_0  v_{p+1}]$.
\end{Lemma}\pagebreak[3]
Next, we will introduce generic annihilator numbers in the exterior algebra.
Let $M \in \M$ and let $v_1,\dots,v_n \in E$ be generic linear forms of $M$.
For $p=1,2,\dots,n$, set
\begin{eqnarray}
A^{(p)}(M)= \big( (v_1,\dots,v_{p-1})M :_M v_{p} \big)\big/ (v_1,\dots,v_{p})M \label{SRM3}
\end{eqnarray}\pagebreak[3]
and
$$\alpha_{p,k}(M)= \dim_K \big(A^{(p)}(M)_k\big).$$
Note that $A^{(p)}(M)= \Ker (\gamma_{0,p-1})$ for $p=2,3,\dots,n$.
These numbers $\alpha_{p,k}(M)$ are constant for a generic choice of linear forms $v_1,\dots,v_n \in E_1$,
and will be called \textit{exterior generic annihilator numbers of $M$}.
In the rest of this section,
we will give the formula to compute the graded Betti numbers of generic initial ideals in the exterior algebra from exterior generic
annihilator numbers.\pagebreak[3]

A monomial ideal $J \subset E$ is said to be \textit{strongly stable}
if $e_S \in J$ and $j \in S$ implies that $e_{(S \setminus \{j\}) \cup \{i\}} \in J$ 
for all $i <j$ with $i \not \in S$.
It is known that generic initial ideals are strongly stable
(\cite[Proposition 1.7]{AHH}).\pagebreak[3]
\begin{Lemma} \label{strongly}
Let $J \subset E$ be a graded ideal.
Then one has
$$\alpha_{p,k}(E/J)= \big|\big\{ e_S \in G\big(\Gin(J)\big)_{k+1}: \max{(S)}=n-p+1\big\}\big|\ \ 
 \mbox{ for }p=1,2,\dots,n,$$
where $|A|$ denotes the cardinality of a finite set $A$ and $G\big(\Gin(J)\big)_{k+1}$ is the set of minimal monomial generators
of $\Gin(J)$ of degree $k+1$.
\end{Lemma}\pagebreak[3]
\begin{proof}
By a generic change of coordinates,
we may assume that $\init (J)=\Gin(J)$ and $v_1,v_2,\dots,v_{p+1}=e_n,e_{n-1},\dots,e_{n-p}$.
Then, by (\ref{SRM3}), we have
\begin{eqnarray*}
A^{(p+1)}(E/J) = \big((e_n,\dots,e_{n-p+1}) + J:_E e_{n-p}\big) \big/ \big((e_n,\dots,e_{n-p}) + J \big),
\end{eqnarray*}
where $p=0,1,\dots,n-1$.\pagebreak[3]
Set
\begin{eqnarray*}
B^{(p+1)}(E/J) = \big((e_n,\dots,e_{n-p+1}) + \init(J):_E e_{n-p}\big) \big/ \big( (e_n,\dots,e_{n-p}) + \init (J) \big).
\end{eqnarray*}\pagebreak[3]
Since we consider the reverse lexicographic order induced by $e_1>\cdots>e_{n}$,
it follows from \cite[Proposition 5.1]{Almost} that
$$
\init\big((e_n,\dots,e_{n-p+1}) +J:_E e_{n-p}\big)= \big((e_n,\dots,e_{n-p+1}) + \init (J):_E e_{n-p}\big)
$$
and
$$
\init \big((e_n,\dots,e_{n-p}) + J\big)= (e_n,\dots,e_{n-p}) + \init(J).
$$
Since $\big((e_n,\dots,e_{n-p+1}) +J:_E e_{n-p}\big) \supset (e_n,\dots,e_{n-p}) + J$
and taking initial ideals does not change Hilbert functions,
it follows that $B^{(p+1)}(E/J)$ and $A^{(p+1)}(E/J)$ have the same Hilbert function.
Thus we have $\alpha_{p,k}(E/J)= \dim_K B^{(p)}(E/J)_k$ for all $k \geq 0$.\pagebreak[3]

Then, to prove the claim,
it is enough to show that the set of monomials
\begin{eqnarray}
&&\big\{[e_S] \in E/ \big( (e_n,\dots,e_{n-p})+ \init (J) \big): \max (S) < n-p,\ e_S \wedge e_{n-p} \in G(\init(J))_{k+1}\big\}
 \label{SRM4}
\end{eqnarray}
forms a $K$-basis of $B^{(p+1)}(E/J)_k$.\pagebreak[3]

If $e_S$ satisfies the condition of (\ref{SRM4}),
then we have $e_S \not \in (e_n,\dots,e_{n-p})+ \init (J)$.
Thus the set (\ref{SRM4}) is indeed the set of $K$-linearly independent monomials belonging to $B^{(p+1)}(E/J)$.
Hence we need to prove that any nonzero monomial $e_S \in B^{(p+1)}(E/J)$ of degree $k$ is contained 
in the set (\ref{SRM4}).

Let $[e_S] \in B^{(p+1)}(E/J)\setminus \{0\}$ be a monomial of degree $k$.
Then we have $e_S \wedge e_{n-p} \in (e_n,\dots,e_{n-p+1})+ \init(J)$.
Also, since $[e_S]$ is not zero,
we have $e_S \not \in (e_n,\dots,e_{n-p})$.
Thus we have $\max(S) <n-p$ and $e_S \wedge e_{n-p} \in \init(J)$.
Since $\init(J)= \Gin(J)$ is strongly stable and $e_S \not \in \init(J)$,
any monomial $e_T \in E$ of degree $k$ which divides $e_S \wedge e_{n-p}$ does not belongs to $\init(J)$.
Thus we have $e_S \wedge e_{n-p} \in G(\init(J))$,
and $[e_S]$ is contained in the set (\ref{SRM4}).
\end{proof}\pagebreak[3]\vspace{-0.7ex}
For a monomial $e_S \in E$, let $m(e_S)= \max (S)$.
If $J\subset E$ is a strongly stable ideal,
then it follows from \cite[Corollary 3.3]{AHH} that
$$ \beta_{ii+k}^E (E/J)= \sum_{p=k+1}^{n} \sum_{  {e_S \in G(J)_{k+1}} \atop{m(e_S) =p}} { p-1+i-1 \choose i-1}
\quad \mbox{ for all } i \geq 1 \mbox{ and all }k \geq 0.
$$
Since every generic initial ideal is strongly stable,
the above equality together with Lemma \ref{strongly} imply
the next lemma.
\begin{Lemma} \label{ginbetti}
Let $J$ be a graded ideal in $E$.
Then one has
$$\beta_{ii+k}^E  \big(E/\Gin(J)\big)= \sum_{p=1}^{n-k} {n-p+i-1 \choose i-1}\alpha_{p,k} (E/ J)$$
for all $i \geq 1$ and all $k \geq 0$.
\end{Lemma}
\section{Rigidity of resolutions over an exterior algebra}
In this section,
we will prove similar results studied in \S 2 for generic initial ideals in an exterior algebra.

Let $M \in \M$.
Throughout this section,
let $v_1,\dots,v_n \in E_1$ be generic liner forms
and write
$H_i(p)_k$, $h_{i,k}(p)$ and $\alpha_{p,k}$ for $H_i(v_1,\dots,v_p;M)_k$, $\dim_K \big(H_i(v_1,\dots,v_p;M)_k\big)$
and $\alpha_{p,k}(M)$ respectively.
Set $\delta_{i,p,k}= \dim_K \big( \Image(\gamma_{i,p})_k\big)$ for $i >0$
and $\delta_{0,p,k}=0$ for all $p,k$.\pagebreak[3]

For  an integer $j \geq 0$,
Lemma \ref{long} yields
the following exact sequence
\begin{eqnarray*}
\cdots
\stackrel{\gamma_{i,p}}{\longrightarrow}
H_i(p)_j
\stackrel{\eta_{i,p}}{\longrightarrow}
H_i(p+1)_j
\stackrel{\psi_{i,p}}{\longrightarrow}
H_{i-1}(p+1)_{j-1}
\stackrel{\gamma_{i-1,p}}{\longrightarrow}
H_{i-1}(p)_j
\longrightarrow
\cdots
\end{eqnarray*}
where $p=1,2,\dots,n-1$.\pagebreak[3]
Then, in the same way as \S 1, we have
\begin{eqnarray}
h_{1,k}(p+1)
&=& h_{1,k}(p)+ \alpha_{p+1,k-1} - \delta_{1,p,k} \label{LRM5}
\end{eqnarray}\pagebreak[3]
and, for $i >1$, we have
\begin{eqnarray}
h_{i,i+k}(p+1) 
&=& h_{i,i+k}(p) + h_{i-1,i+k-1}(p+1) - \{ \delta_{i,p,i+k} + \delta_{i-1,p,i+k}\}. \label{LRM4}
\end{eqnarray}\pagebreak[3]
\begin{Proposition} \label{formula}
With the same notation as above,
one has
\begin{eqnarray}
\quad \ \ h_{i,i+k}(p)&=& \sum_{j=1}^{p} {p-j+i-1 \choose i-1} \alpha_{j,k} \label{LRM6} \\
\nonumber &&-\sum_{s=1}^i \sum_{j=1}^{p-1}
{ p-1 -j +i-1 -(s-1) \choose  i-1 -(s-1)} \bigg\{ \delta_{s,j,s+k} + \delta_{s-1,j,s+k}\bigg\}. 
\end{eqnarray}
\end{Proposition}
\begin{proof}
The proof is quite similar to the proof of Proposition \ref{graded}.
So we will skip some detail calculations.

We use induction on $p$ and $i$.
First, we will show the case $p=1$.
Recall that $C_{\bullet}(v_1;M)$ is the complex
$$\cdots \longrightarrow
C_{i+1}(v_1;M) \stackrel{\partial}{\longrightarrow}
C_{i}(v_1;M)\stackrel{\partial}{\longrightarrow}
C_{i-1}(v_1;M) {\longrightarrow} \cdots
$$\pagebreak[3]
with the differential $\partial(x_1^{(i)})=v_1 x_1^{(i-1)}$.
Thus we have $$H_i(1)_{i+k} \cong \big((M:_M v_1)/v_1 M\big)_k = A^{(1)}(M)_k,$$
and therefore we have $h_{i,i+k}(1)= \alpha_{1,k}$ for all $i \geq 1$ and all $k \geq 0$.
This is equal to the formula (\ref{LRM6}).\pagebreak[3]

Second, we will consider the case $i=1$.
Since we already proved $h_{1,1+k}(1)=\alpha_{1,k}$,
the equation (\ref{LRM5}) says that
\begin{eqnarray*}
h_{1,1+k}(p) 
&=& \{ \alpha_{1,k}+ \cdots + \alpha_{p,k}\}
-\{ \delta_{1,1,1+k} + \cdots + \delta_{1,p-1,1+k}\}
\end{eqnarray*}
which is equal to the formula (\ref{LRM6}).\pagebreak[3]

Finally,
the formula (\ref{LRM6}) for $i>1$ and $p>1$
follows from the equation (\ref{LRM4}) together with the induction hypothesis
in the same way as Proposition \ref{graded}.
\end{proof}

Next, we will show
the following vanishing property of $\Image(\gamma_{i,p})$,
which is an analogue of Lemma \ref{gcor}.

\begin{Lemma} \label{rigid}
Let $i \geq 1$ be a positive integer.
If $\delta_{i,p,k}=0$ for all $1 \leq p \leq n-1$,
then one has $\delta_{i+t,p,k+t}=0$ for all $1 \leq p \leq n-1$ and all $t \geq 0$.
\end{Lemma}

\begin{proof}
It is enough to prove the claim for $t=1$.
Remark that $\delta_{i,p,k}=0$ if and only if
the map $\eta_{i,p}:H_{i}(p)_k \to H_{i}(p+1)_k$ is injective.
Let $\partial_\ell^{(p)}: H_{i+1}(p)_{k+1} \to H_i(p)_k$ be the map defined by
$$\partial_\ell^{(p)}\big([g_0+g_1x_\ell+ g_2x_\ell^{(2)}+\cdots +g_{i+1}x_\ell^{(i+1)}]\big)= 
[g_1+ g_2x_\ell+\cdots +g_{i+1}x_\ell^{(i)}],$$
where $1 \leq \ell \leq p$ and each $g_t$ does not contain the variable $x_\ell^{(s)}$ for all $s \geq 1$.
Thus $\partial_p^{(p)}$ is equal to the map $\psi_{i+1,p-1}$ which appears in Lemma \ref{long}.
Set $\partial^{(p)}= \bigoplus_{\ell=1}^p \partial_\ell^{(p)}$.
Then we have the following commutative diagram.\pagebreak[3]
$$\xymatrix{
  {H_{i+1}(p)_{k+1}}\ar[rr]^{\partial^{(p)}} \ar[d]_{f_p}  &  &  {\bigoplus_{k=1}^p H_i(p)_k}\ar[d]^{h_p}\\
  {H_{i+1}(p+1)_{k+1}} \ar[rr]^{\partial^{(p+1)}} &  &  {\bigoplus_{k=1}^{p+1}H_{i}(p+1)_k} 
}$$\pagebreak[3]
where $h_p$ is the map defined by $h_p(z_1,\dots,z_p)=\big(\eta_{i,p}(z_1),\dots,\eta_{i,p}(z_p),0\big)$
and $f_p$ is the map defined by $f_p(z)= \eta_{i+1,p}(z)$.\pagebreak[3]

Then $\partial^{(1)}$ is injective since $\partial^{(1)}\big([g_{i+1} x_1^{(i+1)}]\big)=[g_{i+1} x_1^{(i)}]$.
Also, by the assumption,
the map $\eta_{i,p}:H_i(p)_k \to H_{i+1}(p+1)_k$ is injective for all $1 \leq p \leq n-1$.
Thus
$h_p$ is injective for all $1 \leq p \leq n-1$.
We will show that if $\partial^{(p)}$ is injective
then $\partial^{(p+1)}$ is also injective.\pagebreak[3]

Set $u \in \Ker ( \partial^{(p+1)})$.
Then we have $\partial^{(p+1)}_{p+1}(u)= \psi_{i+1,p}(u)=0$.
Thus, by the long exact sequence in Lemma \ref{long}, there exists $w \in H_{i+1}(p)$ such that we have
$\eta_{i+1,p}(w)=f_p(w)=u$.
Since $h_p \circ \partial^{(p)}(w)= \partial^{(p+1)} \circ f_p (w)=0$
and $h_p \circ \partial^{(p)}$ is injective by the induction hypothesis,
it follows that $w=0$ and $\partial^{(p+1)}$ is injective.\pagebreak[3]

Now, we proved that
$\partial^{(p)}$ is injective for all $1 \leq p \leq n-1$.
Thus $h_p \circ \partial^{(p)}$ is injective for all $1 \leq p \leq n-1$.
This fact together with the commutative diagram imply that the map 
$\eta_{i+1,p}: H_{i+1}(p)_{k+1} \to H_{i+1}(p+1)_{k+1}$
is injective for all $1 \leq p \leq n-1$.
Hence
we have $\delta_{i+1,p,k+1}=\dim_K\big( \Image (\gamma_{i+1,p})_{k+1}\big) =0$ for all $1 \leq p \leq n-1$.
\end{proof}\pagebreak[3]

Proposition \ref{formula} and Lemma \ref{rigid} imply
the next theorem.
\begin{Theorem} \label{extmodule}
Let $M \in \M$.
Suppose that for some $i>1$ and $k \geq 0$, we have
$\beta^E_{ii+k}(M)=\sum_{j=1}^{n} {n-j+i-1 \choose i-1}\alpha_{j,k}(M)$.
Then
$$\beta^E_{qq+k}(M) = \sum_{j=1}^{n} {n-j+i-1 \choose i-1} \alpha_{j,k}(M)  \ \ \
\mbox{ for all }q \geq 1.$$
\end{Theorem}
\begin{proof}
Since all binomial coefficients in the formula (\ref{LRM6}) are nonzero,
the assumption says that $\delta_{s,p,s+k}=0$ and $\delta_{s-1,p,s+k}=0$ for all $1 \leq s \leq i$ and
 all $1 \leq p \leq n-1$.
Then Lemma \ref{rigid} says that
$\delta_{s,p,s+k}=0$ and $\delta_{s-1,p,s+k}=0$ for all $s \geq 1$ and all $1 \leq p \leq n-1$.
Thus, the statement follow from the formula (\ref{LRM6}).
\end{proof}\pagebreak[3]

Next we consider the case $M=E/J$.
Lemma \ref{strongly} says that, for any graded ideal $J$ of $E$,
one has $\alpha_{j,k}(E/J)=0$ for $j>n-k$.
Thus for any $i \geq 1$ and $k \geq 0$ we have $\sum_{j=1}^n {n-j+i-1 \choose i-1}\alpha_{j,k}(E/J)=\sum_{j=1}^{n-k} {n-j +i-1 \choose i-1}\alpha_{j,k}(E/J)$.
Then the following corollaries follows from Lemma \ref{ginbetti} and Theorem \ref{extmodule} in the same way as in \S 2.
\begin{Corollary} \label{degthm}
Let $J \subset E$ be a graded ideal.
If $\beta^E_{ii+k}(E/J)= \beta^E_{ii+k}\big(E/\Gin(J)\big)$
for some $i >1$ and $k \geq 0$,
then
$$\beta^E_{qq+k}(E/J) =  \beta^E_{qq+k}\big(E/\Gin(J)\big)
\mbox{ for all }q \geq 1.$$
\end{Corollary}\pagebreak[3]
\begin{Corollary} \label{omoituki}
Let $J\subset E$ be a graded ideal.
Then, for a given integer $k$, the graded Betti numbers
 $\beta^E_{ii+k}(E/J)= \beta^E_{ii+k}\big(E/\Gin(J)\big)$ for all $i \geq 1$
if and only if 
$\beta^E_{1,k+1}(E/J)=\beta^E_{1,k+1}\big(E/\Gin(J)\big)$ and $\beta^E_{1,k+2}(E/J)=\beta^E_{1,k+2}
\big(E/\Gin(J)\big)$.
\end{Corollary}\pagebreak[3]
\begin{Remark} {\em Notice that the above Corollary \ref{degthm} and
Corollary \ref{rigidgin} in $\S 2$ are similar.
But as we see Corollary \ref{degthm} is relatively more stronger. We give here an example to show that  in
the case of a
polynomial ring one cannot have the stronger result as in Corollary \ref{degthm}. Consider the ideal
$I=( x_1x_4^2,x_2^3,x_2^2x_3) \subset S =\CC[x_1,x_2,x_3,x_4]$. The minimal graded free resolution of
 $S/I$ and
$S/\Gin(I)$ are given by :}\pagebreak[3]
\begin{eqnarray*}
&&0 \longrightarrow S(-7) \longrightarrow
S(-4) \oplus S^2(-6) \longrightarrow
S^3(-3) \longrightarrow
S \longrightarrow  S/I \longrightarrow
0,
\end{eqnarray*}
and 
\begin{eqnarray*}
&&
0 \longrightarrow S(-7) \longrightarrow
S^2(-4)\oplus S(-5) \oplus S^2(-6) \longrightarrow \\
&&\;\;\;\;\;\;\;\;\;\;\;\;\;\;\;\;\;\;\;\;\;\;\;\;\;\;\; S^3(-3)\oplus S(-4)\oplus S(-5) \longrightarrow S
\longrightarrow S/\Gin I \longrightarrow
0.
\end{eqnarray*}\pagebreak[3]
{\em From above resolutions, we see that
$\beta^S_{2,2+4}(S/I)=\beta^S_{2,2+4}\big(S/\Gin(I)\big)=2$ and ofcourse then
$\beta^S_{3,3+4}(S/I)=\beta^S_{3,3+4}\big(S/\Gin(I)\big)=1$. But the graded Betti number
$\beta^S_{1,1+4}(S/I)=0 \neq 1 = \beta^S_{1,1+4}\big(S/\Gin(I)\big)$.}
\end{Remark}\pagebreak[3]
In the case of exterior algebra, the notions of regularity,
linear resolutions and componentwise linear ideals
are defined in the same way as in the case of polynomial ring.
In \cite[Theorem 2.1]{Stable} it was proved that a graded ideal $J$ in $E$
is componentwise linear if and only if $J$ and $\Gin(J)$ have the same graded Betti numbers. 
Theorem \ref{degthm} and Corollary \ref{omoituki} provide the following  new characterization of componentwise linear ideals
in the exterior algebra.
(See also \cite{NRV} for other characterizations of componentwise linear ideals.)\pagebreak[3]
\begin{Theorem} \label{cw}
A graded ideal $J$ in the exterior algebra $E$ is componentwise linear if and only if
$\beta_{i}^E(E/J) = \beta_i^E \big(E / \Gin(J)\big)$ for some $i \geq 1$.\pagebreak[3]
\end{Theorem}

\begin{proof}
Since $\beta_{ii+k}^E(E/J) \leq \beta_{ii+k}^E \big(E / \Gin(J)\big)$ for all $i\geq 1$ and $k \geq 0$, the 
equality
$\beta_{i}^E(E/J) = \beta_i^E \big(E / \Gin(J)\big)$ implies $\beta_{ii+k}^E(E/J) 
= \beta_{ii+k}^E \big(E / \Gin(J)\big)$
for all $k \geq 0$.
Then Theorem \ref{degthm} and Corollary \ref{omoituki} say that
$\beta_i^E(E/J)=\beta_i^E\big(E/\Gin(J)\big)$ for some $i \ge 1$ if and only if
$J$ and $\Gin(J)$ have the same graded Betti numbers.
Hence the claim follows.
\end{proof}\pagebreak[3]
\section{Linear components and graded Betti numbers}

Throughout this section,
we assume that $R$ is either the polynomial ring $S$ over the field $K$ with char$K=0$
or the exterior algebra $E$ over an infinite field.\pagebreak[3]

First, we will extend Corollaries \ref{rigidgin} and  \ref{degthm} to lexsegment ideals and generic initial ideals
with respect to any term
order.
For a strongly stable ideal $I$ in $R$ and for integers $q =1,\dots,n$ and $k \geq 0$, let
$$m_{\leq q} (I,k)= \big|\big\{u \in I: u \mbox{ is a monomial with } m(u) \leq q \mbox{ and } \deg(u)=k\big\}\big|.$$
\begin{Lemma} \label{trans}
Let $I\subset R$ be a graded ideal and $I' \subset R$ a strongly stable ideal with the same Hilbert function as $I$.
Assume that $I'$ satisfies $m_{\leq q}(I',d) \leq m_{\leq q}\big(\Gin(I),d\big)$ for all $q,d$ and
$\beta_{ii+k}^R(R/I) = \beta_{ii+k}^R (R/{I'})$ for some $i > 1$ and $k \geq 0$.
\begin{itemize}
\item[(i)]
If $R=S$, then one has $\beta_{qq+k}^S(S/I) = \beta_{qq+k}^S (S/{I'})$ for all $q \geq i$.
\item[(ii)]
If $R=E$,
then one has $\beta_{qq+k}^E(E/I) = \beta_{qq+k}^E (E/{I'})$ for all $q \geq 1$.
\end{itemize}
\end{Lemma}\pagebreak[3]

\begin{proof}
We will show the case $R=S$.
(The proof for the case $R=E$ is same.)
It follows from \cite[Proposition 2.3]{B} that,
for any strongly stable ideal $J \subset S$, we have\pagebreak[3]
\begin{eqnarray}
\label{FLC} \beta^S_{ii+j}(S/J) &=&
\dim_K J_{j+1} {n-1 \choose i}\\
\nonumber &&-\sum_{q=i}^{n-1} m_{\leq q}(J,j+1) { k-1 \choose i-1} - \sum_{q=i+1}^{n}m_{\leq q}(J,j){k-1 \choose i}
\end{eqnarray}
for all $i$ and $j$.\pagebreak[3]
(A similar formula for graded Betti numbers over the exterior algebra appears in \cite[Theorem 4.4]{AHH}.)
Then  by (\ref{FLC}) and the assumption, we have
$\beta^S_{ij}(S/I) \leq \beta^S_{ij}\big(S/\Gin(I)\big) \leq \beta^S_{ij}(S/I')$
for all $i,j$.\pagebreak[3]
Thus, by Corollary \ref{rigidgin}, what we must prove is $\beta_{qq+k}^S\big(S/\Gin(I)\big)=
\beta_{qq+k}^S(S/I')$ for all $q \geq i$.
However (\ref{FLC}) and the assumption imply that
$m_{\leq q}\big(\Gin(I),k+1\big)= m_{\leq q}(I',k+1)$ for all $q \geq i$
and $m_{\leq q}\big(\Gin(I),k\big)= m_{\leq q}(I',k)$ for all $q \geq i+1$.
Hence for all $q \geq i$, we have $\beta_{qq+k}^S\big(S/\Gin(I)\big) = \beta_{qq+k}^S (S/{I'})$ 
as desired.\pagebreak[3]
\end{proof}

Let $I \subset R$ be a graded ideal.
We write $\Lex(I)\subset R$ for the unique lexsegment ideal of $R$ with the same Hilbert function as $I$ defined 
in \cite{B} (or \cite{SQlex} for the exterior case)
and $\Gin_\sigma (I)$ for the generic initial ideal of $I$ with respect to a term order $\sigma$.
It is known that $\Lex(I)$ and $\Gin_\sigma(I)$ satisfy the assumption of Lemma \ref{trans}
(see \cite[\S 5]{C} and \cite[\S 5]{NRV}).
Thus we have\pagebreak[3]
\begin{Theorem}
Let $I\subset R$ be a graded ideal, $\sigma$ a term order
and let $J$ be either $\Gin_\sigma(I)$ or $\Lex(I)$.
Suppose that $\beta_{ii+k}^R(R/I) = \beta_{ii+k}^R (R/{J})$ for some $i > 1$.
\begin{itemize}
\item[(i)]
If $R=S$, then one has $\beta_{qq+k}^S(S/I) = \beta_{qq+k}^S (S/{J})$ for all $q \geq i$.
\item[(ii)]
If $R=E$,
then one has $\beta_{qq+k}^E(E/I) = \beta_{qq+k}^E (E/{J})$ for all $q \geq 1$.
\end{itemize}
\end{Theorem}\pagebreak[3]

Next,
we consider when a graded ideal $J$ satisfies
$\beta_{ii+d}^E(E/J)= \beta_{ii+d}^E\big(E/\Gin(J)\big)$ for all $i \geq 1$,
where we fix an integer $d\geq 0$.
The next  lemma follows from \cite{BS} and \cite[Theorem 5.3]{Almost}.\pagebreak[3]
\begin{Lemma} \label{regular}
Let $I\subset R$ be a graded ideal.
Then, $I$ has a linear resolution if and only if $\Gin(I)$ has a linear resolution.
\end{Lemma}\pagebreak[3]

We also require the following.
\begin{Lemma}[Crystallization Principle] \label{green}
Let $I\subset R$ be a graded ideal.
If $I$ is generated by elements of degree $\leq d$
and $\beta^R_{1d+1}\big(R/\Gin(I)\big)=0$,
then $\reg(I) \leq d$.
\end{Lemma}\pagebreak[3]

The Crystallization Principle
was proved by Green \cite[Corollary 2.28]{G}
for generic initial ideals over a polynomial ring,
however, this fact can also be proved for generic initial ideals over an exterior 
algebra in the same way.\pagebreak[3]
\begin{Proposition} \label{linear}
Let $I \subset R$ be a graded ideal.
The following conditions are equivalent.
\begin{itemize}
\item[(i)] $I_{\langle k \rangle}$ has a linear resolution;
\item[(ii)] $\beta^R_{1k+1}(R/I) = \beta_{1k+1}^R\big(R/\Gin(I)\big)$,
that is, the number of elements of degree $k+1$ belonging to the set of minimal generators of $I$ is equal to that of $\Gin(I)$.
\end{itemize}
\end{Proposition}
\begin{proof}
%
Let $\mm$ be the maximal ideal of $R$.
Since $\beta_{1k+1}^R(R/I)$ is the numbers of generators in $G(I)$ of degree $k+1$,
we have
\begin{eqnarray*}
\beta_{1k+1}^R(R/I) &=& \dim_K I_{k+1} - \dim_K (\mideal I_{\langle k \rangle})_{k+1} \\
&=& \dim_K I_{k+1} - \dim_K ( I_{\langle k \rangle})_{k+1}
\end{eqnarray*}
and
\begin{eqnarray*}
\beta_{1k+1}^R(R/\Gin(I)) 
&=& \dim_K (\Gin(I)_{k+1}) - \dim_K \big(\mideal \Gin(I_{\langle k \rangle})\big)_{k+1}.
\end{eqnarray*}\pagebreak[3]
Then, from above equations we have $\beta_{1k+1}^R(R/I) = \beta_{1k+1}^R\big(R/\Gin(I)\big)$ if and only if
$\dim_K (I_{\langle k \rangle})_{k+1}=\dim_K \big(\mm \Gin(I_{\langle k \rangle})\big)_{k+1}.$

Suppose $I_{\langle k \rangle}$ has a linear resolution.
Then, by Lemma \ref{regular} $\Gin(I_{\langle k \rangle})$ has a linear resolution.
Hence
$\dim_K \big( \mm \Gin(I_{\langle k \rangle})\big)_{k+1}
= \dim_K \big(\Gin(I_{\langle k \rangle})\big)_{k+1}
= \dim_K (I_{\langle k \rangle})_{k+1}.$
Hence we have $\beta_{1k+1}^R(R/I) = \beta_{1k+1}^R\big(R/\Gin(I)\big)$ as required. On the other hand,
if $\beta^R_{1k+1}(R/I) = \beta^R_{1k+1}\big(R/ \Gin(I)\big)$,
then $\dim_K \big( \mm \Gin( I_{\langle k \rangle})\big)_{k+1}= 
\dim_K \big(\Gin (I_{\langle k \rangle})\big)_{k+1}$.
This implies
$\beta^R_{1k+1} \big(R/ \Gin (I_{\langle k \rangle})\big)=0$.
Then the Crystallization Principle says that $I_{\langle k \rangle}$ has a linear resolution.
\end{proof}\pagebreak[3]

Now, Theorem \ref{dlinear} immediately follows from the above proposition together with Corollaries \ref{generator} and \ref{omoituki}.
Indeed,
(i) $\Leftrightarrow$ (iii) of Theorem \ref{dlinear} follows from  Corollaries \ref{generator} and \ref{omoituki}.
Also,
(ii) $\Leftrightarrow$ (iii) of Theorem \ref{dlinear} follows from Proposition \ref{linear}.\pagebreak[3]
\begin{Example}
Let $I=(x_1^2,x_2^2,x_1x_2x_3^2,x_3^5)\subset S=\CC[x_1,x_2,x_3]$.
Then we have $$\Gin(I)=
(x_1^2,x_1x_2,x_2^3,x_2^2x_3^2,x_1x_3^4,x_2x_3^5,x_3^6).$$
Then Proposition \ref{linear} says that
$I_{\langle k \rangle}$ has a linear resolution for $k=3,4,7,8,9,\dots.$
In particular,
for $k=4,8,9,10,\dots,$
we have $\beta^S_{ii+k}(I)= \beta^S_{ii+k}\big(\Gin(I)\big)$ for all $i \geq 0$.
\end{Example}\pagebreak[3]
\section{The Cancellation Principle}

Let $K$ be a field of characteristic $0$.
In this section,
we will study the relation between our results in \S 1 and the
Cancellation Principle for generic initial ideals,
which was considered in \cite{G}.
This observation would help us to understand
why we require the assumption $i > 1$ in Corollary \ref{rigidgin}
and why we need to consider  $I_{\langle k \rangle}$ and
$I_{\langle k+1 \rangle}$ in Theorem \ref{dlinear}.\pagebreak[3]

First, we recall what is the Cancellation Principle.

\begin{Lemma}[{\cite[Corollary 1.21]{G}}]
\label{cancell}
Let $I$ be a graded ideal in $S$ and $\sigma$ a term order.
The minimal free resolution of $I$ is obtained from that of
$\init_\sigma (I)$ by cancelling adjacent terms,
in other words,
there exists integers $\tau_{i,i+k}$ with $1 \leq i \leq n-1$ and $k \geq
0$ such that
$$\beta_{ii+k}^S\big(\init_\sigma (I)\big)= \beta_{ii+k}^S (I) + \tau_{i,i+k} +
\tau_{i+1,i+k} \ \ \ \mbox{ for all }i\geq 0 \mbox{ and all } k \geq 0,$$
where we let $\tau_{0,k}=0$ for all $k \geq 0$.
\end{Lemma}\pagebreak[3]

We refer the reader to \cite[Example 1.35]{G}
for further information about the Cancellation Principle.
\medskip

Let $I$ be a graded ideal in $S$.
Then Lemma \ref{cancell} says that
there exists integers $c_{i,i+k}(I)$ with $1 \leq i\leq n-1$ and with $k \geq
0$ such that
$$\beta_{ii+k}^S\big(\gin (I)\big)= \beta_{ii+k}^S (I) + c_{i,i+k}(I) +
c_{i+1,i+k}(I) \ \ \ \mbox{ for all }i\geq 0 \mbox{ and all } k \geq 0,$$
where we let $c_{0,k}(I)=0$ for all $k \geq 0$.
It can be easily verified that
the integers $c_{i,i+k}(I)$ are uniquely determined for a given ideal $I$.
We will call the integer $c_{i,i+k}(I)$ the
\textit{$(i,i+k)$-th  cancellation number of $I$}.\pagebreak[3]

\begin{Example}
Let $I=(x_1^3,x_1^2x_2,x_1x_2^2,x_2^3,x_1^2x_3,x_1x_3x_4)\subset S=\CC[x_1,x_2,x_3,x_4]$.
Then we have $\gin(I)=(x_1^3,x_1^2x_2,x_1x_2^2,x_2^3,x_1^2x_3,x_1x_2x_3,x_1x_3^3)$.
The minimal free resolution of $I$ is
\[
0 \longrightarrow
S(-5) \oplus S(-6) \longrightarrow
S^6(-4) \oplus S(-5) \longrightarrow
S^6(-3) \longrightarrow
I \longrightarrow
0,
\]
and that of $\gin(I)$ is
\[
0 \to
S^2(-5) \oplus S(-6) \to
S^7(-4) \oplus S^2(-5) \to
S^6(-3) \oplus S(-4) \to
\gin(I) \to
0.
\]
Hence we have $c_{1,4}(I)=1$, $c_{2,5}(I)=1$ and
all other  cancellation numbers of $I$ are $0$.
\end{Example}\pagebreak[3]

In \S 2,
we already proved that (see Proposition \ref{graded} and (\ref{TT}))
\begin{eqnarray*}
\beta^S_{ii+k}(I)=\beta^S_{ii+k}\big(\gin(I)\big)-
\! \sum_{(a,b)\in A_{i+1,n}}\left[\binom{n-b-1}{i-a+1}\delta_{a,b,a+k-1}
+\binom{n-b-1}{i-a}\delta_{a,b,a+k}\right],
\end{eqnarray*}
where $\delta_{a,b,a+k}=\dim_K\big( \Image \varphi_{a,b}\big)_{a+k}$
and where $\varphi_{a,b}$ is the map which appears in the long exact sequence (\ref{Kos}).
This formula enables us to
write the cancellation numbers
in terms of the Koszul homology of generic linear forms.\pagebreak[3]
\begin{Lemma}
\label{can}
With the same notation as above,
one has
$$c_{i,i+k}(I)=\sum_{(a,b)\in A_{i+1,n}}\binom{n-b-1}{i-a}\dim_K\big( \Image
\varphi_{a,b}\big)_{a+k} \ \ \ \mbox{ for all }i\geq 0 \mbox{ and all } k \geq
0.$$
\end{Lemma}
\begin{proof}
For all $i \geq 0$ and all $k \geq 0$, we set $C_{i,i+k}=\sum_{(a,b)\in A_{i+1,n}}\binom{n-b-1}{i-a}
\dim_K\big(\Image
\varphi_{a,b}\big)_{a+k}$ and $C'_{i,i+k}=\sum_{(a,b)\in A_{i+1,n}}\binom{n-b-1}{i-a+1}\dim_K\big(\Image
\varphi_{a,b}\big)_{a+k-1}$.
Then we have $$\beta^S_{ii+k}(I)=\beta^S_{ii+k}\big(\Gin(I)\big)-C_{i,i+k}-C'_{i,i+k}.$$
Notice that we only need to show that $C'_{i,i+k}=C_{i+1,i+k}.$
Recall that, in the proof of Theorem \ref{jaan},
we already proved that
$$A_{i+2,n} \setminus \big\{(i+2,b):b \leq n-1\big\}
=A_{i+1,n} \setminus \big\{(i-n+b+1,b):b\leq n-1\big\}.$$\pagebreak[3]
Now,
since the binomial $\binom{n-b-1}{i-a+1}=0$ for all $(a,b) \in \big\{(i+2,b):b \leq n-1\big\}$
and for all $(a,b) \in\big\{(i-n+b+1,b):b\leq n-1\big\}$, we have\pagebreak[3]
\begin{eqnarray*}
C_{i+1,i+k}&=&\sum_{(a,b)\in A_{i+2,n}}\binom{n-b-1}{i-a+1}\dim_K\big( \Image
\varphi_{a,b}\big)_{a+k-1}\\
&=&\sum_{(a,b)\in A_{i+2,n}\setminus \{(i+2,b):b \leq n-1\}}
\binom{n-b-1}{i-a+1}\dim_K\big(\Image \varphi_{a,b}\big)_{a+k-1}\\
&=&\sum_{(a,b)\in A_{i+1,n} \setminus \{(i-n+b+1,b):b\leq n-1\}}
\binom{n-b-1}{i-a+1}\dim_K\big(\Image \varphi_{a,b}\big)_{a+k-1}\\
&=&\sum_{(a,b)\in A_{i+1,n} }
\binom{n-b-1}{i-a+1}\dim_K\big(\Image \varphi_{a,b}\big)_{a+k-1}\\
&=& C'_{i,i+k}.
\end{eqnarray*}
This concludes the proof.
\end{proof}\pagebreak[3]
By using Lemma \ref{can},
we can prove an analogue of Corollaries \ref{rigidgin} and \ref{generator}.\pagebreak[3]

\begin{Theorem} \label{crigid}
Let $I$ be a graded ideal in $S$.
If $c_{i,i+k}(I)=0$ for some $i \geq 1$ and  $k \geq 0$,
then one has $c_{q,q+k}(I)=0$ for all $q \geq i$.
\end{Theorem}
\begin{proof}
It suffices to show the case $q=i+1$.
Remark that $\dim_K \big(\Image \varphi_{a,b}\big)_{a+k}=0$
if and only if $\big(\mm H_a(b)\big)_{a+k}=0$.
In the proof of Theorem \ref{jaan},
we proved that  if $\dim_K\big( \Image \varphi_{a,b}\big)_{a+k}=0$
for all $(a,b) \in A_{i+1,n} \setminus \big\{(i+1,b): b \leq n-1\big\}$,
then $\dim_K\big(\Image \varphi
_{a,b}\big)_{a+k}=0$ for all
$(a,b) \in A_{i+2,n} \setminus \big\{(i+2,b):b \leq n-1\big\}$.
Then, since ${n-b-1 \choose i-a+1} =0$ for any
$(a,b) \in  \big\{ (i+2,b):b \leq n-1 \big\}$,
Lemma \ref{can} says that $c_{i+1,i+1+k}(I)=0$.
\end{proof}\pagebreak[3]
\begin{Corollary}\label{clinear}
Let $I$ be a graded ideal in $S$.
Then $c_{i,i+k}(I)=0$ for all $i\geq 1$
if and only if $I_{ \langle k \rangle }$ has a linear resolution.
\end{Corollary}\pagebreak[3]
\begin{proof}
Since the graded Betti number
$\beta_{ 0,k+1 }^S\big(\gin(I)\big)= \beta_{ 0,k+1 }^S(I)+
c_{1,1+k}(I)$,
we have $\beta_{0,k+1}^S\big(\gin(I)\big)= \beta_{0,k+1}^S(I)$ if and only if
$c_{1,1+k}(I)=0$.
However, by Theorem \ref{crigid},
we have $c_{1,1+k}(I)=0$ if and only if $c_{i,i+k}(I)=0$ for all $i \geq 1$.
Also, by Proposition \ref{linear},
we have $\beta_{0,k+1}^S\big(\gin(I)\big)= \beta_{0,k+1}^S(I)$ if and only if
$I_{\langle k \rangle}$ has a linear resolution.
Thus the assertion follows.
\end{proof}\pagebreak[3]

Observe that Theorems \ref{crigid} and Corollary \ref{clinear}
are stronger than Corollaries \ref{rigidgin} and \ref{generator}.
Indeed,
Corollary \ref{rigidgin}
immediately follows from Theorem \ref{crigid},
since the graded Betti numbers $\beta^S_{ii+k}(I)=\beta^S_{ii+k}\big(\gin(I)\big)$
if and only if $c_{i,i+k}(I)=0$ and $c_{i+1,i+k}(I)$=0.\pagebreak[3]

We also remark the next fact which follows from
Lemma \ref{can}.\pagebreak[3]
\begin{Corollary}
Let $I$ be a graded ideal in $S$.
Assume that $I_{\langle k \rangle}$ has a linear resolution.
\begin{itemize}
\item[(i)]
If $\beta^S_{q,q+k+2}(I)= \beta^S_{q,q+k+2}\big(\gin(I)\big)$,
then $\beta^S_{q+1,q+k+2}(I)= \beta^S_{q+1,q+k+2}\big(\gin(I)\big)$;
\item[(ii)] If $\beta^S_{q,q+k-1}(I)=\beta^S_{q,q+k-1}\big(\gin(I)\big)$,
then $\beta^S_{q-1,q+k-1}(I)=\beta^S_{q-1,q+k-1}\big(\gin(I)\big)$.
\end{itemize}
\end{Corollary}\pagebreak[3]
\begin{proof}
By Corollary \ref{clinear},
we have $c_{\ell, \ell+k}(I)=0$ for all integers $\ell \geq 1$.
Then, we have the graded Betti numbers 
$\beta^S_{q+1,q+k+2}\big(\gin(I)\big)=\beta^S_{q+1,q+k+2}(I)+ c_{q+1,q+k+2}(I)$
and
$\beta^S_{q-1,q+k-1}\big(\gin(I)\big)=\beta^S_{q-1,q+k-1}(I)+ c_{q,q+k-1}(I)$.
On the other hand,
if the graded Betti number $\beta^S_{q,q+k+2}(I)= \beta^S_{q,q+k+2}\big(\gin(I)\big)$
then we have $c_{q+1,q+k+2}(I)=0$.
Also,
if $\beta^S_{q,q+k-1}(I)=\beta^S_{q,q+k-1}\big(\gin(I)\big)$
then we have $c_{q,q+k-1}(I)=0$.
Thus the assertion follows.
\end{proof}\pagebreak[3]

As for any graded ideal $I$,
$I_{\langle 1 \rangle}$ always has a linear resolution,
it follows that if $\beta_{q,q+3}^S \big(\gin(I)\big)=\beta^S_{q,q+3}(I)$ then
we have $\beta_{q+1,q+3}^S \big(\gin(I)\big)=\beta^S_{q+1,q+3}(I)$.\pagebreak[3]

Since it is not difficult  to find  the Betti numbers of a strongly stable ideal $J$,
one may expect to find all possible Betti numbers of graded ideals  $I$
 such that $\Gin(I)=J$ by using  Betti numbers of $J$ and by considering
all possible cancellations.
However,
this problem is far reaching as pointed out in \cite[Example 1.35]{G}.\pagebreak[3]
\bigskip

\noindent
\textbf{Thanks:}
All of the examples that we have presented in the paper
are computed by the computer algebra system \cocoa \ \cite{coco}.
We also mention that computations of generic initial ideals are
done by a random choice of matrices.

\end{document}